\documentclass[10pt]{article}

\usepackage[T1]{fontenc}
\usepackage[latin1]{inputenc}
\usepackage{mathtools,amsthm,amssymb}

\usepackage{mathrsfs}

\usepackage{authblk}

\usepackage{hyperref}
\usepackage[margin=1in]{geometry}

\newtheorem{theo}{Theorem}[section]
\newtheorem{cor}[theo]{Corollary}
\newtheorem{lem}[theo]{Lemma}
\newtheorem{prop}[theo]{Proposition}

\numberwithin{equation}{section}

\def \RR{\mathbb{R}}
\def \SS{\mathbb{S}}
\def \CC{\mathbb{C}}
\newcommand{\Sa}{ {\cal S }}
\newcommand{\Ga}{ {\cal G }}
\newcommand{\Ma}{ {\cal M }}

\def\tr{\mbox{\rm Tr}}

\newcommand{\vertiii}[1]{{\left\vert\kern-0.25ex\left\vert\kern-0.25ex\left\vert #1
    \right\vert\kern-0.25ex\right\vert\kern-0.25ex\right\vert}}

\begin{document}

\title{An Explicit Floquet-Type Representation of Riccati Aperiodic Exponential Semigroups}

\author[$1$]{Adrian N. Bishop}
\author[$2$]{Pierre Del Moral}
\affil[$1$]{{\small University of Technology Sydney (UTS); and Data61 (CSIRO), Australia}}
\affil[$2$]{{\small INRIA, Bordeaux Research Center, France; and UNSW Sydney, Australia}}
\date{}

\maketitle

\begin{abstract}
The article presents a rather surprising Floquet-type representation of time-varying transition matrices associated with a class of nonlinear matrix differential Riccati equations. The main difference with conventional Floquet theory comes from the fact that the underlying flow of the solution matrix is aperiodic. The monodromy matrix associated with this Floquet representation coincides with the exponential (fundamental) matrix associated with the stabilizing fixed point of the Riccati equation. The second part of this article is dedicated to the application of this representation to the stability of matrix differential Riccati equations. We provide refined global and local contraction inequalities for the Riccati exponential semigroup that depend linearly on the spectral norm of the initial condition. These refinements improve upon existing results and are a direct consequence of the Floquet-type representation, yielding what seems to be the first results of this type for this class of models.
\end{abstract}

\section{Introduction}\label{sec-models}

The aim of this article is to provide the analog of a Floquet normal form, commonly used to represent the fundamental matrices of periodic linear dynamical systems, in the case of matrix Riccati differential semigroups. These new exponential semigroup formulae shed new light on the regularity properties of matrix differential Riccati equations. Moreover, we present this Floquet normal form as an interesting and explicit example of a case in which the exponential semigroup operator of a non-periodic differential equation can be expressed explicitly.

\subsection{Notation}

We introduce some standard notation. We denote by $\Ma_{r}=\RR^{r\times r}$ the set of $(r\times r)$-square matrices with real entries and $r\geq 1$. We let $\Sa_r\subset \Ma_{r}$ denote the subset of symmetric matrices, $\Sa_r^0\subset\Sa_r$ the subset of positive semi-definite matrices, and $\Sa_r^+\subset \Sa_r^0$ the subset of positive definite matrices. Given $B\in \Sa_r^0-\Sa_r^+$ we denote by $B^{1/2}$ a (non-unique) but symmetric square root of $B$ (given by a Cholesky decomposition). When $B\in\Sa_r^+$ we always choose the principal (unique) symmetric square root. We write $A^{\prime}$ to denote the transposition of a matrix $A$, and $A_{sym}=(A+A^{\prime})/2$ to denote the symmetric part of $A\in\Ma_{r}$. 

We equip the set $\Ma_{r}$ with the spectral norm $\Vert A \Vert=\sqrt{\lambda_{max}(AA^{\prime})}$ where $\lambda_{max}(\cdot)$ denotes the maximal eigenvalue. The minimal eigenvalue is denoted by $\lambda_{min }(\cdot)$. Let $\tr(A)=\sum_{1\leq i\leq r}A(i,i)$ denote the trace operator. We also denote by $\mu(A)=\lambda_{max}(A_{sym})$ the logarithmic norm and by $\varsigma(A):=\max_i{\left\{\mbox{\rm Re}[\lambda_i(A)]\right\}}$ the spectral abscissa. Note that $\varsigma(A)\leq\mu(A)$.

\subsection{Matrix Differential Riccati Equations}

Given some given matrices $(A,R,S)\in (\Ma_{r}\times\Sa^0_r\times \Sa^0_r)$ we denote by $\phi_{s,t}(P_s)=P_t$ the (forward in time) flow associated with the matrix differential Riccati equation,
 \begin{align}
	\partial_tP_t~=~\Lambda(P_t) ~:=&~ AP_t + P_tA^{\prime}+R-P_tSP_t \nonumber\\
	:=&~ (A-P_tS)P_t+P_t(A-P_tS)^{\prime}+ R+P_tSP_t,~\qquad P_0 = Q  \label{def-Riccati-drift} 
\end{align}	
with $0\leq s\leq t$ and the convention $\phi_{t}(Q):=\phi_{0,t}(Q)$. We suppose $Q\in\Sa^0_r$, the pair $(A,S^{1/2})$ observable and the pair $(A,R^{1/2})$ is controllable.

Let $E_{s,t}(Q)$ denote the exponential semigroup associated with the matrix flow $t\mapsto (A-\phi_t(Q)S)$; that is
the solution for any $0\leq s\leq t$ of the equations
\begin{equation}\label{syst-E-def}
 \partial_t E_{s,t}(Q)=(A-\phi_t(Q)S)\,E_{s,t}(Q)\quad\mbox{\rm and}\quad
\partial_s E_{s,t}(Q)=-E_{s,t}(Q)\,(A-\phi_s(Q)S)
\end{equation}
with $E_{s,s}(Q)=I$ and where we often write $E_{t}(Q)$ for $E_{0,t}(Q)$. The observability and controllability conditions ensure \cite{Lancaster1995} the existence and uniqueness of a matrix $P_{\infty}\in\Sa_r^+$ with
\begin{equation}\label{CARE}\Lambda(P_{\infty})=0\quad\mbox{\rm and such that}\quad
B:=A-P_{\infty}S\quad\mbox{\rm is a  stable matrix (a.k.a. Hurwitz matrix).}
\end{equation}	
So $\varsigma(B)<0$. Observe that for any initial matrix $Q\in\Sa^0_+$ and any time horizons $s,t\in\RR_+$ we have
\begin{equation}\label{Est}
E_{s,s+t}(Q)=E_{t}(\phi_s(Q))\quad \mbox{\rm and}\quad
E_{s,s+t}(P_{\infty})=E_{t}(P_{\infty})=e^{tB}
\end{equation}	
The implicit solution for $\phi_t(Q)$ is given by
\begin{align}
	\phi_t(Q) \,=\, E_{t}(Q)\,Q\,E_{t}'(Q) + \int_{0}^t\,E_{s,t}(Q)\left[R + \phi_{s}(Q)\,S\,\phi_{s}(Q)\right]E_{s,t}(Q)'\,ds
\end{align}
and it follows that the stability of this solution is intimately connected to the stability of the exponential semigroup (\ref{syst-E-def}). Nevertheless, most of the Lipschitz and contraction-type estimates developed in the literature are derived using Lyapunov and Gronwall-type techniques involving rather crude constants. A review on this subject in the context of filtering can be found in~\cite{ap-2016}. In this work, we provide an explicit Floquet normal form for $E_{t}(Q)$ which is surprising since the flow $\phi_t(Q)$ is aperiodic. We then explain how this Floquet-type representation can be used to refine those contraction properties and results on $E_{t}(Q)$.

\subsection{Some History and Convergence Literature}

Bucy \cite{bucy2} originally studied a number of global properties of the matrix differential Riccati equation (\ref{def-Riccati-drift}). In particular, he proved that solutions exist for all time when the initial condition is positive semi-definite, he proved a number of important monotonicity properties, along with bounds\footnote{The original upper and lower bounds on the (time-varying) Riccati equation given in \cite{bucy2} were particularly elegant in appearance; being given in terms of the relevant observability and controllability Gramians. However, as noted in \cite{hitz72}, there was a crucial (yet commonly made) error in the proof which invalidated the result as given. This error was repeated (and/or overlooked) in numerous subsequent works. A correction \cite{bucy72corrected} was noted in a reply to \cite{hitz72}; see Bucy's reply \cite{bucy72remarks} and a separate reply by Kalman \cite{kalman72}. However, a complete reworking of the result did not appear in entirety, it seems, until much later in \cite{delyon2001}.} on the solution stated in terms of the controllability and observability Gramians. Bucy \cite{bucy2} studied when the solution of the autonomous Riccati equation converges to a solution of an associated (fixed-point) algebraic Riccati equation, and finally he proved exponential stability of the time-varying Kalman-Bucy filter along with an exponential forgetting property of the associated Riccati equation.

Associated with the differential Riccati equation is the (fixed-point) algebraic equation (\ref{CARE}) whose solution(s) correspond to the equilibrium point(s) of the corresponding differential equation. This algebraic equation was studied by Bucy in \cite{bucy72} and it was shown that there exists an unstable negative definite solution (in addition to the desired positive definite equilibrium). A detailed study of the algebraic Riccati equation was given by J.C. Willems \cite{willems71} who considered characterising all solutions. Bucy \cite{bucy75} later considered the so-called structural stability of these solutions. Detectability and stabilisability conditions are necessary and sufficient for a unique stabilising positive semi-definite solution of the algebraic equation \cite{kucera72}. See also \cite{kucera72,callier81,wimmer85} and the early review paper \cite{kucera73} for related literature. A (marginally) stable solution of the algebraic Riccati equation exists under detectability conditions; see \cite{Molinari73,Molinari77,wimmer85,Poubelle86} and \cite{VanHandel2009}. The discussion in \cite[Chapter 2 and 3]{Bittanti91} is also of general interest here; as is \cite{Lancaster1995}.

Returning to the differential Riccati equation (\ref{def-Riccati-drift}), convergence to a stabilising fixed-point was studied extensively in \cite{callier81}, where generalised convergence conditions were given. Detectability and stabilisability conditions are sufficient with a suitable initial condition; e.g. see also the seminal text \cite{Kwakernaak72}. We also note the early paper \cite{wonham68} that studied convergence and dealt further with a generalized version of the Riccati equation with a linear perturbation term. A geometric analysis of the differential Riccati equation and its solution(s) is given in \cite{shayman86}. In \cite{gevers85} sufficient conditions are given such that the solution of the differential Riccati equation at any instant is stabilising; see also \cite{Poubelle88}. In \cite{nicolao92,callier95,Park97} convergence to a (marginally) stable solution was studied under relaxed conditions (with necessity also explored). In \cite{VanHandel2009}, given only detectability, the stability of a time-varying ``closed loop'', e.g. $(A-\phi_t(Q)S)$, is proven even when the limiting Riccati solution, or $(A-P_\infty S)$, is only marginally stable.

Finally, we point to \cite{ap-2016} for a number of convergence results and contraction estimates on the (time-varying) Riccati flow $\phi_t(Q)$ and on the exponential semigroup $E_{t}(Q)$. As noted, we seek to refine many of these convergence estimates given the explicit Floquet-type normal representation of $E_{t}(Q)$ presented later.

\subsection{Solution Techniques and Forms}

In order to estimate the stability properties of $\phi_t(Q)$ or $E_{t}(Q)$ it would be useful to represent these quantities in some rather explicit form. Appropriately chosen representations also aid in robustness and sensitivity analysis (e.g. in terms of perturbations of the system model $(A,R,S)$). The literature on representations of $\phi_t(Q)$ is too broad to adequately survey but in this subsection we highlight briefly some work in this direction with further pointers to the literature. See the early reviews in, e.g., \cite{anderson-moore,kenney} for initial guidance. Representations of $E_{t}(Q)$ are much less studied, and the Floquet-type normal form presented herein is one of the main contributions of this work. 

From the numerical viewpoint, matrix differential Riccati equations (\ref{def-Riccati-drift}) may be solved for the flow $\phi_t(Q)$ using direct Runge-Kutta numerical integration \cite{kenney,Dieci92}. Ensemble Kalman-Bucy filtering \cite{evensen03} provides a mean field interacting particle approximation of the Riccati flow $\phi_t(Q)$ in terms of sample covariance matrices associated with a Kalman-type differential equation. This probabilistic approach for numerical evaluation of the flow is also rather general, and the resulting sample covariance evolution is known to satisfy (\ref{def-Riccati-drift}) with an additional (additive) diffusion term that can be made as small as desired at the expense of more computation (i.e. more particles). The stability analysis of these matrix diffusion processes rely on the stability properties of (\ref{def-Riccati-drift}); see for instance~\cite{Bishop/DelMoral/Bruno2012,apa-2017,BishopDelMoralMatricRicc,delmoral16a}.

More explicit formulae for the flow $\phi_t(Q)$ itself are numerous in the literature. The Bernoulli substitution method \cite{reid72,bucy2,Davison73,kenney,mceneaney}, and Chandrasekhar decompositions \cite{Kailath72,Lindquist74,Lainiotis76,kenney}, lead to several related expressions for the flow $\phi_t(Q)$. For example, the Bernoulli substitution approach has numerous variants \cite{kenney}, but the broad idea is to express the solution of the Riccati equation as the ratio of two matrices (e.g. associated with the spectral decomposition of a Hamiltonian matrix). A system of linear differential equations on an augmented state space (double the original system size) is then solved for the individual matrices forming the ratio. This strategy is rather general, but there is no intuitive way to relate with some precision, the ratio-solution nor the spectrum of the Hamiltonian to the parameters of the model. The article~\cite{ferrante} provides one description of this couple of matrices in terms of the steady state matrix $P_{\infty}$. 

Numerous other expressions for $\phi_t(Q)$ exist; e.g. methods based on Lyapunov transformations are also popular, see \cite{anderson-moore,Leipnik85,Nazarzadeh98,Nguyen2010}. The literature here is too broad to adequately survey. See other representations for $\phi_t(Q)$ in, e.g., \cite{Rusnak88,Bittanti91,jodar1991closed,mceneaney,ferrante}. The literature referenced in these articles points to further forms of $\phi_t(Q)$ (or variations of more common forms discussed here). 

In \cite{prach2015infinite,prach-thesis} an expression for the flow $\phi_t(Q)$ is given in which explicit connections with the system model and the steady state solution can be made. This expression for $\phi_t(Q)$ is derived herein as an immediate consequence of the Floquet-type representation on the exponential semigroup $E_{t}(Q)$ given later.

Coming to the exponential semigroup $E_{t}(Q)$, we find that explicit expressions for $E_{t}(Q)$ are less studied. The definition (\ref{syst-E-def}) implies that $E_{t}(Q)$ is also the exponential semigroup (i.e. the fundamental matrix) associated with the time-varying linear system, $\dot{x}_t \,=\, (A-P_tS)x_t$. 

Interestingly, the contractive properties $E_{t}(Q)$ imply the stability of both the preceding linear (time-varying) differential system and the time-invariant matrix differential Riccati equation (\ref{def-Riccati-drift}). The study of time-varying linear differential equations is difficult and intuition does not generally carry over from the time-invariant case. For example, in contrast with time-invariant systems, the eigenvalue placement in general (arbitrary) time-varying coefficient matrices does not yield information on the stability of the system. Well known counterexamples of unstable systems associated to time-varying matrices with negative eigenvalues are given in~\cite{coppel1978stability}. Conversely, stable systems defined by time-varying matrices with positive eigenvalues are given in~\cite{wu1974note}. Most stability results for systems like $\dot{x}_t \,=\, (A-P_tS)x_t$, in general, rely on the design of judicious Lyapunov functions, often described by ad-hoc time varying matrix inequalities \cite{khalil}. This approach is also difficult in practice. See alternatively \cite{rosenbrock63,Desoer69,coppel1978stability,Amato93,solo94,Kumar,Ilchmann,Mullhaupt07} for special cases.

One idea for simplifying the stability analysis of time-varying linear systems is to study explicit representations of the exponential semigroup. However, this approach is typically not available since expressions for the fundamental matrix are generally not obtainable. However, we may contrast this with the time-invariant case where the fundamental matrix is just a standard matrix exponential \cite{coppel1978stability}, and more generally with periodic systems (as discussed and exploited subsequently). 

From a mathematical viewpoint, it is tempting to integrate sequentially the differential equations (\ref{syst-E-def}), to obtain an explicit description of $E_{s,t}(Q)$ in terms of Peano-Baker series~\cite{Peano,baker1905}, see also~\cite{Brockett,Frazer,Ince}. Another natural strategy to express the semigroup as a true
 matrix exponential 
\begin{equation*}
	\textstyle E_{s,t}(Q)=\exp{\left[\oint_s^t \, (A-\phi_u(Q)S)\, du\right]}
\end{equation*}
where $\oint_s^t \cdot\, du$ denotes the traditional Magnus series expansion involving iterated integrals on the Lie algebra generated by the matrices $(A-\phi_u(Q)S)$, with $s\leq u\leq t$. For more details on these exponential expansions we refer to~\cite{blanes,magnus}. In practical terms, the use of Peano-Baker and/or exponential Magnus series in the study of the stability properties of time-varying linear dynamical systems is rather limited.

\subsection{Floquet Theory}

It is well known that the fundamental matrix (i.e. the exponential semigroup) of time-varying linear systems {\em with periodic matrix coefficients} can be represented in terms of the exponential of a matrix (a.k.a. the monodromy matrix) up to an invertible and bounded matrix flow. This feature is the main subject of Floquet theory~\cite{floquet}.  
  
The central idea is to design a natural change of coordinates that transforms any time varying but periodic systems to traditional linear and time homogeneous models. In this situation, the stability of the system is characterized in terms of the spectrum of the exponential monodromy matrix associated to the system; see e.g.~\cite{brown,dacunha,gorcek} and the references therein.
This elegant and powerful approach provides sharp, as well as necessary and sufficient conditions, for the stability of periodic linear systems systems. 
 
Besides the fact that the flow of Ricatti matrices $t\in\RR_+:=[0,\infty[\mapsto A_t:=A-P_tS$ is not periodic and the matrices $A_sA_t\not=A_tA_s $ don't commute ($s\neq t$) we have the following rather surprising theorem.

\begin{theo}[Floquet-type Representation]
For any time horizon $t\geq 0$ and any $Q\in\Sa_r^0$ we have Riccati exponential  semigroup formula
\begin{equation}\label{floquet-type-ric}
E_{t}(Q)= e^{t B }~\CC_t(Q)^{-1}\quad\mbox{with}\quad
\sup_{t\geq 0}{\Vert \CC_t(Q)^{-1}\Vert}<\infty\quad\mbox{and}\quad \CC_t(P_{\infty})=I
\end{equation}
where $B$ is defined in (\ref{CARE}) and $\CC_t:\Sa^0_r\mapsto \Ga l_r$ is some flow of linear and invertible functionals, where
$\Ga l_r\subset \Ma_r$ denotes the general linear group of invertible matrices. Clearly, $E_{t}(P_\infty)=e^{t B }$.
\end{theo}
A more precise description of the matrix functional $\CC_t$, including some uniform estimates is provided in Section~\ref{main-models-1}; see Theorem~\ref{second-theo-intro}. This result is interesting as besides the fact that the evolution $\phi_t(Q)$ and the matrix flow $\CC_t(Q)$ are not periodic, most of the stability properties discussed in Floquet theory still apply to this class of model. Section~\ref{main-models-1} also contains an explicit form of the Riccati flow $\phi_t(Q)$ derived easily from the Floquet representation of $E_{t}(Q)$, and some explicit expressions of the limiting matrix associated with some particular models $(A,R,S)$.

Note that a Floquet representation of the linear system (on an augmented state space, double the original system size) associated with the Bernoulli substitution method \cite{bucy2,Davison73,kenney} for expressing $\phi_t(Q)$ has been explored in \cite{abou-kandil03}. This study holds for matrix differential Riccati equation with periodic coefficients \cite{Bittanti91,abou-kandil03} (of which the time-invariant case here is of course a special case). However, we emphasise that the Floquet representation on the level of the linear system associated with the Bernoulli substitution does not yield a Floquet form on the level of the Riccati flow $\phi_t(Q)$ itself; i.e. it does not represent a Floquet form of $E_{t}(Q)$.

In Section~\ref{main-models-2}, we illustrate the impact of the Floquet-type formula (\ref{floquet-type-ric}) in the stability analysis of the matrix differential equation (\ref{def-Riccati-drift}). We provide a series of refined spectral and Lipschitz 
contraction estimates on $E_{t}(Q)$. To get some intuition here on the power or usefulness of this representation, we recall that $\phi_t$ is a smooth matrix functional with a first order Fr\'echet derivative defined for any $Q\in\Sa_r^0$ and $H\in\Sa_r$ by
$$
(\ref{floquet-type-ric}) ~~~~\Longrightarrow~~~~
\nabla\phi_t(Q)\cdot H\,=\,E_{t}(Q)\,H\,E_{t}(Q)^{\prime}\,=\,e^{tB}\,\CC_t(Q)^{-1}\,H\,(\CC_t(Q)^{\prime})^{-1}e^{tB^{\prime}}
$$
Higher order Fr\'echet derivatives can also be found in~\cite{apa-2017}. In addition, for any $Q_1,Q_2\in\Sa^0_r$ we have the decomposition
$$
\phi_t(Q_1)-\phi_t(Q_2)=E_{t}(Q_1)\,(Q_1-Q_2)\,E_{t}(Q_2)^{\prime}
$$
The above formula is a direct consequence of the polarisation formula (\ref{polarization-formulae}). These decompositions combined with the Floquet representation (\ref{floquet-type-ric}) allows one to develop a refined analysis of the local and global contraction properties of the Riccati semigroup functionals $Q\mapsto\phi_t(Q)$. To the best of our knowledge, the uniform and quadratic-type estimates presented here are the first results of this type for this class of models.

\section{A Floquet Exponential Semigroup Formula}\label{main-models-1}

Consider the Gramian flow and a related steady state matrix,
\begin{eqnarray}
\SS_t&:=&\int_0^t\, e^{sB^{\prime}}
S~e^{sB}\,ds ~~~\longrightarrow_{t\rightarrow\infty}~~~ \SS_{\infty}>0\quad\mbox{\rm and set}\quad P^{-}_{\infty}:=P_{\infty}-\SS_{\infty}^{-1}
\label{def-SS-OO}
\end{eqnarray}
The function $\SS:t\in\RR_+\mapsto \SS_t\in \Sa^+_r$ is positive definite for any $t>0$ due to the observability assumption \cite{prach2015infinite} and it is increasing w.r.t. the Loewner order. The matrix $P^{-}_{\infty}<0 $ is the unique negative-definite steady state of the Riccati equation (cf. Proposition~\ref{prop-negative-steady-state} given later). Note that the existence of the negative-definite fixed point $P^{-}_{\infty}$ follows from the observability and controllability assumption; see \cite[Chapter 3]{Bittanti91} and \cite{Lancaster1995}.

Consider now the linear matrix functional
\begin{equation}\label{inter-f-0}
\CC_t~:~Q\in \Sa_r^0\mapsto
\CC_t(Q):=\left[(\SS_t^{-1}-\SS^{-1}_{\infty})+(Q-P_{\infty}^-)\right]~\SS_t~\in \Ga l_r
\end{equation}
Note also that rearranging gives,
\begin{equation}\label{inter-f-1}
(\ref{def-SS-OO})~~~\Longrightarrow~~~
\CC_t(Q)=I+(Q-P_{\infty})\,\SS_t~~~\Longrightarrow~~~ \CC_t(P_{\infty})=I
\end{equation}
The invertibility of $\CC_t$ is immediate. 

Our main result takes the following form.

\begin{theo}\label{second-theo-intro}
The  state transition matrix semigroup $E_t(Q)$ is given for any $t>0$ in closed form by the exponential formulae (\ref{floquet-type-ric}).

In addition, for any $\delta>0$ we have the uniform estimates
\begin{eqnarray}
\sup_{t\geq \delta}{\Vert \CC_t(Q)^{-1}\Vert} \,\leq \chi_{\delta}&\Longrightarrow&~~ \forall t\geq \delta\quad
\Vert E_t(Q)\Vert \,\leq~\chi_{\delta}\,\Vert E_t(P_{\infty})\Vert\label{estimate-E-noQ}\\
\sup_{t\geq 0}{\Vert \CC_t(Q)^{-1}\Vert}\,\leq \chi(Q)&\Longrightarrow&~~ \forall t\geq 0\quad
\Vert E_t(Q)\Vert \,\leq~\chi(Q)\,\Vert E_t(P_{\infty})\Vert\label{estimate-EQ}
\end{eqnarray}
where again we highlight that $E_t(P_{\infty})=e^{tB}$ and where
$$
\chi(Q):= \Vert P_{\infty}^-\Vert^{-1}
\left[\Vert P_{\infty}-P_{\infty}^-\Vert+\Vert Q-P_{\infty}\Vert\right]\quad\mbox{\rm and}\quad
\chi_{\delta}:=\left[{\lambda_{ min}\left(\SS_{\delta}\right)
\lambda_{ min}\left(-P_{\infty}^-\right)}\right]^{-1}
$$
\end{theo}
The proof is provided in Section~\ref{second-theo-proof}. We leave it as open, the problem of establishing an analogue of Theorem \ref{second-theo-intro} under milder detectability and stabilisability conditions \cite{kucera72,Kwakernaak72,callier81,Lancaster1995} on the model $(A,R,S)$.

Using the decomposition
\begin{equation}
\Lambda(Q_1)-\Lambda(Q_2)
\,=\,(A-Q_1S)(Q_1-Q_2)+(Q_1-Q_2)(A-Q_2S)^{\prime}\label{polarization-formulae}
\end{equation}
we readily check the following corollary.
\begin{cor}\label{bernstein}
For any $t>0$ we have the closed form Lipschitz type matrix formula
\begin{equation}\label{eqn-cor-decomp}
\phi_t(Q_1)-\phi_t(Q_2)=E_t(P_{\infty})~\CC_t(Q_1)^{-1}(Q_1-Q_2)~\left(E_t(P_{\infty})~\CC_t(Q_2)^{-1}\right)^{\prime}
\end{equation}
\end{cor}

Applying Corollary~\ref{bernstein} with $Q_2=P_{\infty}$ and using (\ref{inter-f-1}) we recover the Bernstein-Prach-Tekinalp formula~\cite{prach2015infinite,prach-thesis} given by,
$$
\phi_t(Q)=P_{\infty}+ e^{tB}~\left[I+(Q-P_{\infty})~\SS_t\right]^{-1}(Q-P_{\infty})~e^{tB^{\prime}}
$$

\subsection{An Explicit Example Case}

Explicit formulae for the fixed points $(P_{\infty}^-,P_{\infty})$ in terms of $(A,R,S)$ are difficult to obtain\footnote{For example, in~\cite{choiclosed} an explicit formula for $P_{\infty}$ is given under the following particular conditions, $S>0$, and,
$$
R=S^{-1}-AS^{-1}A^{\prime}\quad \mbox{\rm and}\quad 
S^{-1/2}(S^{-1}A^{\prime}-AS^{-1})S(S^{-1}A^{\prime}-AS^{-1})S^{-1/2}+4S^{-2}\geq 0
$$
The article~\cite{choiclosed} is based on the Moser-Veselov interpretation of some class of Riccati equations developed in~\cite{cardoso2003moser}. Related closed form expressions when $R=0$ and with diagonal matrices $A$ are presented in~\cite{rojas-2010}.}. However, when $S>0$ and $SA=A^{\prime}S\geq 0$,
the fixed point matrices $(P_{\infty}^-,P_{\infty})$ are given by the formulae 
\begin{equation}
 P_{\infty}^-\,=\,\left[A-(A^2+RS)^{1/2}\right]S^{-1}<0<P_{\infty}=\left[A+(A^2+RS)^{1/2}\right]S^{-1} \label{formula-example-Pinfty}\end{equation}
which implies $B=-(A^2+RS)^{1/2} = -\frac{1}{2}(P_{\infty}-P_{\infty}^-)S$. The exponential formula (\ref{floquet-type-ric}) resumes to
\begin{equation}\label{formula-example-Pinfty-bis}
	E_t(Q)\,=\,e^{tB}\,(P_{\infty}-P_{\infty}^-)\, \left[e^{2tB}\,(P_{\infty}-P_{\infty}^-)+(Q-P_{\infty}^-)\left(I-e^{2tB}\right)\right]^{-1}
\end{equation}
The derivation of (\ref{formula-example-Pinfty}) and (\ref{formula-example-Pinfty-bis}) are provided in the Appendix \ref{formula-example-Pinfty-ref}.

\section{Applications of the Floquet Form: Refined Contraction Estimates}\label{main-models-2}

Given a function $F:\Sa^0_r\mapsto \Ma_r$ defined on the open set $\Sa^0_r$ we set
\begin{equation*}
\vertiii{F}\,:=\,\sup{\left\{\Vert F(Q)\Vert~:~Q\in  \Sa^0_r\right\}}\quad\mbox{\rm and}\quad \mbox{\rm Lip}(F)\,:=\,\sup{\left\{\mbox{lip}_F(Q_1,Q_2)\,:\,Q_1,Q_2\in \Sa^0_r\right\}} \,\in\, [0,\infty]
\end{equation*}
with the local Lipschitz constant: $
\mbox{\rm lip}_F(Q_1,Q_2):=\inf{\left\{c~:~
\Vert F(Q_1)- F(Q_2)\Vert\leq c~\Vert Q_1-Q_2\Vert\right\}}$.

\subsection{Some Existing Contraction Estimates}

For any $Q\in \Sa_{r}^0$ we have the following results (e.g. see \cite{bucy2,ap-2016}),
\begin{equation}\label{exp-stab-phi}
\forall t\geq \delta >0\quad\Pi_{-,\delta}\leq \phi_{t}(Q)\leq \Pi_{+,\delta}\qquad\mbox{and}\qquad \forall t\geq 0\quad
 \Vert e^{tB}\Vert \,\leq\, \alpha\,e^{-\beta t}
\end{equation}
for some matrices $\Pi_{-,\delta},\Pi_{+,\delta}\in \Sa^+_{r}$ and some $\alpha,\beta>0$ that all depend on the model $(A,R,S)$. The parameters $(\alpha,\beta)$ can be made explicit in terms of the spectrum $B$. For instance, applying Coppel's inequality (cf. Proposition 3 in~\cite{coppel1978stability}), for any $t\geq 0$ we can choose $\alpha=({a}/{\gamma})^{r-1}$ and $\beta=-(1-\gamma)\varsigma(B)$ for any $0<\gamma<1$ with $a:=2{\Vert B\Vert}/{\vert\varsigma(B)\vert}$. We can pick $\alpha=1$ if it holds that $\beta=-\mu(B)>0$. There is also a natural change of basis such that $P_\infty=I$ which also allows one to choose $\alpha=1$.

Using Lyapunov techniques one can show (viz. e.g. \cite{bucy2,ap-2016}) that for any $t\geq 0$, $s\geq\delta>0$ we have
\begin{equation}\label{exp-stab-phi-bis}
\vertiii{ E_{s,s+t}}~\leq~ \alpha_\delta \,e^{-\beta_\delta t} 
\end{equation}
for some $\alpha_\delta,\beta_\delta>0$ that depend also on the model $(A,R,S)$. Note in \cite{ap-2016} the assumption that $R$ is invertible is taken; as is common in filtering applications where $R$ relates to the signal noise covariance matrix. In \cite{ap-2016}, the constants in this estimate are explored quite explicitly with this assumption of $R$ invertible (or more technically with $R$ \emph{or} $S$ invertible). However, qualitatively, this bound (\ref{exp-stab-phi-bis}) easily holds more generally under just controllability/observability of $(A,R,S)$; see \cite{bucy2,Bittanti91}, or even under weaker stabilisability/detectability conditions; see \cite{kucera72,Kwakernaak72}. 

Corollary~4.9 in~\cite{ap-2016} provides for any $s\geq 0$ exponential decays 
\begin{equation}\label{exp-stab-phi-bis-alltimes}
\Vert E_{s,s+t}(Q)\Vert~\leq~ \exp{(\alpha\Vert Q\Vert)}\,\,e^{-\beta t}\quad \mbox{\rm  }
\end{equation}
for some $\alpha,\beta>0$, different to those in (\ref{exp-stab-phi}), and that depend on $(A,R,S)$ with again $R$ (or $S$) invertible.

Note that the contraction inequality (\ref{exp-stab-phi-bis-alltimes}) holds for all time horizons. However, in \cite{ap-2016} the penalty we pay for allowing $s\rightarrow0$ in (\ref{exp-stab-phi-bis-alltimes}) is the exponential constant in $\Vert Q\Vert$. 

In the next subsection, we easily (and quantitatively) refine these estimate as a consequence of the Floquet representation and the refined spectral estimates on $\CC_t(Q)$ in Theorem \ref{second-theo-intro}. We also only ask for the model $(A,R,S)$ to be controllable and observable. 

Note also that it is desirable to relate the decay of $E_{s,s+t}(Q)$ to the decay at the fixed point $e^{tB}$ (since as $t\rightarrow\infty$ it is clear that we cannot do better).

\subsection{Refining the Contraction Estimates}

With Theorem~\ref{second-theo-intro} and Corollary~\ref{bernstein}, the exponential constant $\exp{(\alpha\Vert Q\Vert)}$ in the preceding (existing) contraction estimates can be replaced by $c\,(1+\Vert Q\Vert)$, and without any resort to requiring $R>0$ as used in \cite{ap-2016} when explicitly computing constants. Throughout this subsection, $c<\infty$ is some parameter whose values depends on the matrices $(A,R,S)$ and may change from line to line.

Recall throughout that $E_{t}(P_\infty)=e^{t B }$ and $e^{t B }$ is bounded and decays as in (\ref{exp-stab-phi}). We have the following corollary.

\begin{cor}
For any $s,t\geq 0$ we have
\begin{equation}\label{Est-bound}
  \Vert E_{s,s+t}(Q)\Vert ~\leq~ c\,\left(1+\Vert Q\Vert\right)\,\Vert E_{t}(P_{\infty})\Vert
\end{equation}
In addition, for any $t\geq 0$ and $s\geq \delta>0$ we have the uniform estimate
\begin{equation}\label{Est-bound-bis}
  \vertiii{ E_{s,s+t}} \leq c\,\Vert E_{t}(P_{\infty})\Vert\qquad \mbox{and}\qquad
  \vertiii{ E_{t,t+s}} \leq~\chi_{\delta}\,\Vert E_{s}(P_{\infty})\Vert
\end{equation}
\end{cor}
\proof
Combining (\ref{Est}) with (\ref{estimate-EQ}) and Proposition 4.3 in~\cite{ap-2016} for any $t\geq s\geq 0$ we check that
\begin{equation*}
  \Vert E_{s,t}(Q)\Vert
  \,\leq\, \left( \Vert P_{\infty}^-\Vert^{-1}
\left[\Vert P_{\infty}-P_{\infty}^-\Vert+\Vert \phi_s(Q)-P_{\infty}\Vert\right]\right)\,\Vert E_{t-s}(P_{\infty})\Vert ~~\Longrightarrow~~ (\ref{Est-bound})
 \end{equation*}
 The proof of the l.h.s. in (\ref{Est-bound-bis}) follows from the preceding line and the l.h.s. estimate in (\ref{exp-stab-phi}). Using (\ref{estimate-E-noQ}) for any
 $t-s\geq \delta>0$ also have
$$
E_{s,t}(Q)\,=\,E_{t-s}(\phi_s(Q))\quad\Longrightarrow\quad \Vert E_{s,t}(Q)\Vert\,\leq~\chi_{\delta}\,\Vert E_{t-s}(P_{\infty})\Vert
$$
This ends the proof of the corollary. \qed

This refinement, which replaces the exponential constant $\exp{(\alpha\Vert Q\Vert)}$ in, e.g., (\ref{exp-stab-phi-bis-alltimes}) with $c\,(1+\Vert Q\Vert)$ in (\ref{Est-bound}), may be important in practical analysis. For example, traditional backward semigroup analysis involves interpolating paths of the form $s\mapsto\phi_{s,t}(\phi^\epsilon_s(Q))$ for some approximating flow $\phi^\epsilon_s(Q)$ as $\epsilon\rightarrow0$. For stochastic approximations, the convergence analysis requires estimating moments of the matrix exponentials $E_{t}(\phi^\epsilon_s(Q))$. Using (\ref{exp-stab-phi-bis-alltimes}) to estimate such moments would require the control of the exponential moments of $\phi^\epsilon_s(Q)$, whereas using (\ref{Est-bound}) only requires conventional moments. For example, in ensemble Kalman filtering theory, the former often do not exist while the later do; see e.g. \cite{Bishop/DelMoral/Bruno2012,apa-2017,BishopDelMoralMatricRicc}.

We consider the parameters
$$
\left(\chi_{\phi,\delta},\chi_{\phi}(Q_1,Q_2)\right):=\left(\chi_{\delta}^2,\chi(Q_1)\chi(Q_2)\right)\quad\mbox{\rm and}\quad
\left(\chi_{E,\delta},\chi_{E}(Q_1,Q_2)\right):=\Vert\SS_{\infty}\Vert~\left(\chi_{\phi,\delta},\chi_{\phi}(Q_1,Q_2)\right)
$$
In this notation, combining Corollary~\ref{bernstein} with the decomposition
\begin{eqnarray*}
E_t(Q_1)-E_t(Q_2)&=&E_t(P_{\infty})~\CC_t(Q_1)^{-1}~\left[\CC_t(Q_2)-\CC_t(Q_1)\right]~\CC_t(Q_2)^{-1}\\
&=&E_t(P_{\infty})~\CC_t(Q_1)^{-1}~(Q_2-Q_1)~\SS_t~\CC_t(Q_2)^{-1}
\end{eqnarray*}
we readily check the following corollary.
\begin{cor}
For any $t\geq \delta>0$ we have the uniform contraction inequalities
\begin{equation}\label{cor-Et}
\mbox{\rm Lip}(\phi_t)\leq  \chi_{\phi,\delta}~\Vert E_t(P_{\infty})\Vert^2\quad\mbox{and}\quad \mbox{\rm Lip}(E_t)\leq  \chi_{E,\delta}~
\Vert E_t(P_{\infty})\Vert
\end{equation}
In addition, for any time horizon $t\geq 0$ we also have the local Lipschitz estimates
\begin{equation}\label{cor-Et-2}
\mbox{\rm lip}_{\phi_t}(Q_1,Q_2)\leq \chi_{\phi}(Q_1,Q_2)~\Vert E_t(P_{\infty})\Vert^2\quad\mbox{and}\quad
\mbox{\rm lip}_{E_t}(Q_1,Q_2)\leq  \chi_{E}(Q_1,Q_2)~\Vert E_t(P_{\infty})\Vert~
\end{equation}
\end{cor}

Observe that 
$$
	\chi_{\phi}(Q_1,Q_2)\vee \chi_{E}(Q_1,Q_2)~\leq~ c\,(1+\Vert Q_1\Vert^2+\Vert Q_2\Vert^2)
$$
It follows that the estimates in the preceding corollary improve those presented in Section 4.3 in~\cite{ap-2016} that depend on some exponential Lipschitz constant of the form $ \exp{(\alpha\,(\Vert Q_1\Vert+\Vert Q_2\Vert))}$.

\section{Proof of Theorem~\ref{second-theo-intro}}\label{second-theo-proof}

\begin{prop}\label{prop-negative-steady-state}
The matrix $P^{-}_{\infty}<0 $ is the unique negative-definite steady state of the Riccati equation. In addition, for any $s< t$ we have
$$
0< \SS_s<\SS_t\leq \SS_{\infty}\quad\mbox{and}\quad P_{\infty}<\SS_{\infty}^{-1}
$$
with the positive matrices $\SS_t$, $\SS_{\infty}$ defined in (\ref{def-SS-OO}).
\end{prop}
\proof
The proof of l.h.s. assertion is immediate.
Observe that 
$$
\SS_{\infty}B+B^{\prime}~\SS_{\infty}+S=0
$$
By definition of the steady state $P_{\infty}>0$ we have
$$
\begin{array}{l}
P_{\infty}^{-1}A+A^{\prime}P_{\infty}^{-1}+P_{\infty}^{-1}RP_{\infty}^{-1}-S~=~P_{\infty}^{-1}
B+B^{\prime}P_{\infty}^{-1}+P_{\infty}^{-1}RP_{\infty}^{-1}+S~=~0\\
\\
\qquad\qquad\displaystyle\Longrightarrow\qquad
(P_{\infty}^{-1}-\SS_{\infty})B+B^{\prime}(P_{\infty}^{-1}-\SS_{\infty})+P_{\infty}^{-1}RP_{\infty}^{-1}~=~0\\
\\
\qquad\qquad\displaystyle\Longrightarrow \qquad P_{\infty}^{-1}~=~\SS_{\infty}+\int_0^{\infty}
e^{sB^{\prime}}\,P_{\infty}^{-1}RP_{\infty}^{-1}\,e^{sB}\,ds~>~\SS_{\infty}
\end{array}$$
This implies that $P_{\infty}<\SS_{\infty}^{-1}$. Now, applying polarisation formula
\begin{equation*}
\Lambda(Q_1)-\Lambda(Q_2) ~=~(A-Q_2S)(Q_1-Q_2)+(Q_1-Q_2)(A-Q_2S)^{\prime}-(Q_1-Q_2)S(Q_1-Q_2)
\end{equation*}
to the pair $(Q_1,Q_2)=(P^{-}_{\infty},P_{\infty})$ we have
$$
\begin{array}{l}
B~\SS_{\infty}^{-1}+\SS_{\infty}^{-1}B^{\prime}+\SS_{\infty}^{-1}S\,\SS_{\infty}^{-1}~=~0\\
\\
\Longleftrightarrow~~~~ B\,(P^{-}_{\infty}-P_{\infty})+(P^{-}_{\infty}-P_{\infty})B^{\prime}-(P^{-}_{\infty}-P_{\infty})S\,(P^{-}_{\infty}-P_{\infty})~=~\Lambda(P^{-}_{\infty})~=~0
\end{array}$$
The uniqueness of $P^{-}_{\infty}<0$ is proven in \cite{bucy72}. This completes the proof of the proposition. \qed

Observe that
$$
\CC_t(Q)~~\longrightarrow_{t\rightarrow\infty}~~\CC_{\infty}(Q)\,=\,I+(Q-P_{\infty})
(P_{\infty}-P^-_{\infty})^{-1}\,=\,(Q-P_{\infty}^-)(P_{\infty}-P^-_{\infty})^{-1}
$$
As product of positive-definite matrices may fail to be positive-definite, the mapping $Q\mapsto\CC_t(Q)$ is not necessarily a positive map. Nevertheless, $\CC_t(Q)$ is invertible and we have  
the formula
\begin{eqnarray*}
\CC_t(Q)^{-1}&=&\SS_t^{-1}\left[(\SS_t^{-1}-\SS^{-1}_{\infty})+(Q-P_{\infty}^-)~\right]^{-1}\nonumber\\
&=&I-(Q-P_{\infty})~\left[(\SS_t^{-1}-\SS^{-1}_{\infty})+(Q-P_{\infty}^-)~\right]^{-1} \qquad\left(\Longleftarrow (\ref{inter-f-1})\right)
\end{eqnarray*}
\begin{lem}
For any $\delta>0$ we have the estimates
\begin{eqnarray*}
t\geq \delta&\Longrightarrow&\Vert \CC_t(Q)^{-1}\Vert\leq \left[{\lambda_{ min}\left(\SS_{\delta}\right)
\lambda_{ min}\left((\SS_t^{-1}-\SS^{-1}_{\infty})-P_{\infty}^-\right)}\right]^{-1}\\
t\leq \delta&\Longrightarrow&
\Vert \CC_t(Q)^{-1}\Vert\leq
1+\Vert Q-P_{\infty}\Vert~/~\lambda_{ min}\left(\SS_{\delta}^{-1}-\SS^{-1}_{\infty}-P_{\infty}^-\right)
\end{eqnarray*}
In addition, we have the l.h.s uniform estimate in (\ref{estimate-EQ}).
\end{lem}
\proof
we have
$$
0<\left[
(\SS_t^{-1}-\SS^{-1}_{\infty})+(Q-P_{\infty}^-)~\right]^{-1}\leq \left((\SS_t^{-1}-\SS^{-1}_{\infty})-P_{\infty}^-\right)^{-1}
$$
This implies that
$$
0<\lambda_{ max}\left(\left[
(\SS_t^{-1}-\SS^{-1}_{\infty})+(Q-P_{\infty}^-)~\right]^{-1}
\right)
\leq \frac{1}{\lambda_{ min}\left((\SS_t^{-1}-\SS^{-1}_{\infty})-P_{\infty}^-\right)}
$$
from which we conclude that
$$
\Vert\SS_t^{-1}\Vert\,\Vert\left[
(\SS_t^{-1}-\SS^{-1}_{\infty})+(Q-P_{\infty}^-)~\right]^{-1}\Vert ~\leq~
\frac{1}{\lambda_{ min}\left(\SS_t\right)\lambda_{ min}\left((\SS_t^{-1}-\SS^{-1}_{\infty})-P_{\infty}^-\right)}
$$
On the other hand, we also have
$$
t\geq \delta \qquad\Longrightarrow\qquad
\SS_{\delta} \,\leq\, \SS_t \qquad\Longrightarrow\qquad \lambda_{ min}\left(\SS_t\right)\,\geq\, \lambda_{ min}\left(\SS_{\delta}\right)
$$
This ends the proof of the first estimate which implies the l.h.s in (\ref{estimate-E-noQ}). Observe that
$$
\begin{array}{l}
\displaystyle
\Vert
I-(Q-P_{\infty})~\left[(\SS_t^{-1}-\SS^{-1}_{\infty})+(Q-P_{\infty}^-)~\right]^{-1}\Vert^2\\
\\
\qquad\qquad\qquad\displaystyle\leq ~1+{\Vert Q-P_{\infty}\Vert^2}/{\lambda_{ min}\left((\SS_t^{-1}-\SS^{-1}_{\infty})+(Q-P_{\infty}^-)\right)^2}
\end{array}
$$
On the other hand for any $0\leq t\leq \delta$ we have
$$
\SS_t^{-1}\geq \SS^{-1}_{\delta} \qquad\Longrightarrow\qquad
\lambda_{ min}\left(\SS_t^{-1}-\SS^{-1}_{\infty}-P_{\infty}^-\right) \,\geq \,
\lambda_{ min}\left(\SS_{\delta}^{-1}-\SS^{-1}_{\infty}-P_{\infty}^-\right)
$$
This ends the proof of the second estimate. This shows that for any $\delta>0$  we have
$$
\sup_{t\geq 0}{\Vert \CC_t(Q)^{-1}\Vert} ~\leq~
\left[1+\frac{\Vert Q-P_{\infty}\Vert}{\lambda_{ min}\left(\SS_{\delta}^{-1}-\SS^{-1}-P_{\infty}^-\right)}\right]\vee
\frac{1}{\lambda_{ min}\left(\SS_{\delta}\right)
\lambda_{ min}\left(-P_{\infty}^-\right)}
$$
Letting $\delta$ tend to $\infty$ we conclude that
$$
\sup_{t\geq 0}{\Vert \CC_t(Q)^{-1}\Vert} ~\leq~ \frac{1}{\lambda_{ min}\left(-P_{\infty}^-\right)}~\left(
\left[\lambda_{ min}\left(-P_{\infty}^-\right)+\Vert Q-P_{\infty}\Vert\right]\vee
\frac{1}{\lambda_{ min}\left(\SS_{\infty}\right)}\right)
$$
This implies that
$$
\sup_{t\geq 0}{\Vert \CC_t(Q)^{-1}\Vert} ~\leq~
\left[1+\frac{\Vert Q-P_{\infty}\Vert}{\lambda_{ min}\left(-P_{\infty}^-\right)}\right]\vee
\frac{\lambda_{ max}\left(P_{\infty}-P_{\infty}^-\right)}{\lambda_{ min}\left(-P_{\infty}^-\right)}
$$
which implies the l.h.s. in (\ref{estimate-EQ}). This ends the proof of the lemma. \qed

Now we come to the {\bf Proof Theorem~\ref{second-theo-intro}:}

\noindent We set
$$
F_t(Q)\,=\,E_t(P_{\infty})~\CC_t(Q)^{-1} \quad\mbox{\rm and}\quad P_t\,:=\,P_{\infty}+F_t(Q) (Q-P_{\infty})E_t(P_{\infty})^{\prime}
$$
Observe that
$$
\begin{array}{l}
\partial_t\CC_t(Q)\,=\,(Q-P_{\infty})\,E_t(P_{\infty})^{\prime}S\,E_t(P_{\infty})\\
\\
\qquad\Longrightarrow\qquad E_t(P_{\infty})\,\partial_t\CC_t(Q)^{-1}\,=\,-F_t(Q)(Q-P_{\infty})E_t(P_{\infty})^{\prime}\,S\,F_t(Q)
\end{array}
$$
This implies that
\begin{eqnarray*}
\partial_tF_t(Q)&=&B\,F_t(Q)-F_t(Q)(Q-P_{\infty})E_t(P_{\infty})^{\prime}\,S\,F_t(Q) \\
&=&\left(A-\left[P_{\infty}+F_t(Q)\,(Q-P_{\infty}) \,E_t(P_{\infty})^{\prime}~\right]S\right)\,F_t(Q)\,=\,
(A-P_tS)\,F_t(Q)
\end{eqnarray*}
The polarisation formula in (\ref{polarization-formulae}) implies that
\begin{equation*}
\partial_tP_t\,=\,(A-P_tS)\,(P_t-P_{\infty})+(P_t-P_{\infty})\,B^{\prime}\,=\,\Lambda(P_t)
\end{equation*}
and by the solution uniqueness of the Riccati equation we conclude that
$$
P_t=\phi_t(Q) \quad\Longrightarrow\quad F_t(Q)=E_t(Q)
$$ 
The proof of the theorem is finished. \qed

\appendix
\section{Appendix: Derivation of (\ref{formula-example-Pinfty}) and (\ref{formula-example-Pinfty-bis})}
\label{formula-example-Pinfty-ref}

We set
$$(\overline{A},\overline{R},\overline{S})\,=\,\left(S^{1/2}AS^{-1/2},S^{1/2}RS^{1/2},I\right)\qquad\mbox{\rm  and} \qquad
\overline{P}\,=\,S^{1/2}PS^{1/2}
$$
In this notation, we have
\begin{equation*}
SA\,=\,A^{\prime}S\geq 0 \qquad\Longleftrightarrow\qquad \overline{A}\in \Sa_r^0~~\Longrightarrow ~~ \overline{A}^2\in \Sa_r^0
\end{equation*}
On the other hand whenever $P=P_{\infty}$ we have
\begin{eqnarray*}
\Lambda(P)=0&\Longleftrightarrow&~~
\overline{P}^2-\overline{A}\,\overline{P}-\overline{P}\,\overline{A}-\overline{R}\,=\,\left(\overline{P}-\overline{A}\right)^2-(\overline{A}^2+\overline{R})=0\\
&\Longleftrightarrow&~~\overline{P}\,=\,\overline{A}+(\overline{A}^2+\overline{R})^{1/2}~\left(>\overline{A}\right)\\
&\Longleftrightarrow&~~
A-PS\,=\,-S^{-1/2}\,(\overline{A}^2+\overline{R})^{1/2}S^{-1/2}~\left(<0\right)
\end{eqnarray*}
The last assertion comes from the uniqueness property of the Riccati fixed point $P_{\infty}$ s.t. $A-P_{\infty}S$ is stable. 
We conclude that
$$
P_{\infty}\,=\,AS^{-1}+S^{-1/2}(S^{1/2}[A^2S^{-1}+R]S^{1/2})^{1/2}S^{-1/2}
$$
On the other hand, for any $S,Q>0$  we have 
$$
S^{1/2}QS^{1/2}\,=\,S^{1/2}QSS^{-1/2}\,=\,\left[S^{1/2}(QS)^{1/2}S^{-1/2}\right]^2
$$ 
where $(QS)^{1/2}$ has positive eigenvalues. Then, we have the formula
$$
\left[S^{1/2}QS^{1/2}\right]^{1/2}\,=\,S^{1/2}(QS)^{1/2}S^{-1/2} \qquad\Longrightarrow\qquad
S^{-1/2}\left[S^{1/2}QS^{1/2}\right]^{1/2}S^{-1/2}\,=\,(QS)^{1/2}S^{-1}
$$
from which we conclude that
$$
P_{\infty}\,=\,\left[A+\left(A^2+RS\right)^{1/2}\right]S^{-1}
$$
We set $$
\overline{P}^-_{\infty} \,:= \,S^{1/2}P^-_{\infty}S^{1/2}\qquad \mbox{\rm and}\qquad \overline{P}_{\infty}\,:=\,S^{1/2}P_{\infty}S^{1/2}
$$ 
In this notation, we have
$$
P^{-}_{\infty}\,=\,P_{\infty}-\SS^{-1}_{\infty} \qquad\Longleftrightarrow\qquad \overline{P}^{-}_{\infty}\,=\,\overline{P}_{\infty}-\left(S^{-1/2}~\SS_{\infty}~S^{-1/2}\right)^{-1}
$$
On the other hand, recalling that $\overline{A}<\overline{P}_{\infty}$ we have
$$
\begin{array}{l}
S^{-1/2}\,\SS_{\infty}\,S^{-1/2}\\
\\
\qquad\displaystyle=\int_0^{\infty}\,S^{-1/2}\exp{\left[t\left(S^{1/2}\left[\overline{A}-\overline{P}_{\infty}\right]S^{-1/2}\right)\right]}\,S\,\exp{\left[t\left(S^{-1/2}\left[\overline{A}-\overline{P}_{\infty}\right]S^{1/2}\right)\right]}S^{-1/2}\,dt\\
\\
\qquad\displaystyle=\int_0^{\infty}\,\exp{\left(2t\left[\,\overline{A}-\overline{P}_{\infty}\right]\right)}\,dt~=~2^{-1}\,\left[~\overline{P}_{\infty}-\overline{A}\right]^{-1}
\end{array}
$$
We conclude that
$$
\overline{P}^{-}_{\infty}\,=\,2\overline{A}-\overline{P}_{\infty}\,=\,\overline{A}-(\overline{A}^2+\overline{R})^{1/2}~(<0)
$$
and in other words we have
\begin{eqnarray*}
P_{\infty}^-&=&AS^{-1}-S^{-1/2}\,\left(S^{1/2}\left[A^2S^{-1}+R\right]S^{1/2}\right)^{1/2}S^{-1/2}~=~\left[A-\left(A^2+RS\right)^{1/2}\right]S^{-1}~<~ 0
\end{eqnarray*}
This ends the proof of  (\ref{formula-example-Pinfty}). Arguing as above we also have
$$
\begin{array}{l}
S^{-1/2}\,\SS_{t}\,S^{-1/2}\\
\\
\qquad\displaystyle=\int_0^{t}\,\exp{\left(2s\left[~\overline{A}-\overline{P}_{\infty}\right]\right)}\,ds
~=~-2^{-1}\,\left[\,\overline{A}-\overline{P}_{\infty}\,\right]^{-1}\,\left(I-\exp{\left(2t\left[\,\overline{A}-\overline{P}_{\infty}\right]\right)}\right)
\end{array}
$$
which that implies
$$
\begin{array}{l}
\SS_{t}~
\displaystyle
=-2^{-1}\,S~B^{-1}\,\left(I-e^{2tB}\right)~=~-2^{-1}\,\left(I-e^{2tB}\right)\,B^{-1}\,S
\end{array}
$$
Observe that
$$
P_{\infty}-P_{\infty}^-~=~2\left(A^2+RS\right)^{1/2}S^{-1}~=~-2B S^{-1}
$$
Now (\ref{formula-example-Pinfty-bis}) is a direct consequence of (\ref{floquet-type-ric}). \qed

\end{document}